\documentclass[11pt]{amsart}
\renewcommand{\bar}{\overline}
\usepackage{amsmath,amsthm,amscd,euscript}
\def \r{\mathbb R}

\def \z{\mathbb Z}

\renewcommand{\cot}{\mathop{\rm ctg}\nolimits}

\DeclareMathOperator{\ls}{lS} \DeclareMathOperator{\ld}{ld}
\DeclareMathOperator{\lil}{l\ell}

\newtheorem{theorem}{Theorem}[section]

\newtheorem{statement}[theorem]{Statement}
\newtheorem{proposition}[theorem]{Proposition}
\newtheorem{corollary}[theorem]{Corollary}
\theoremstyle{remark}
\newtheorem{remark}[theorem]{Remark}
\theoremstyle{definition}
\newtheorem{definition}[theorem]{Definition}

\newtheorem{problem}{Problem}
\newtheorem{conjecture}[problem]{Conjecture}

\title{On invariant M\"obius measure and Gauss-Kuzmin face distribution.}
\author{O.~N.~Karpenkov}
\date{29 September 2006}

\thanks{ The work is partially supported by RFBR SS-1972.2003.1
grant, by NWO-RFBR 047.011.2004.026 (RFBR 05-02-89000-NWO\_a)
grant, by RFBR 05-01-02805-CNRSL\_a grant, and by RFBR grant
05-01-01012a.}

\thanks{E-mail: karpenk@mccme.ru.}

\begin{document}
\input epsf

\maketitle

\begin{flushright}
{\it Dedicated to my mentor

Vladimir Igorevich Arnold.}
\end{flushright}


\sloppy \normalsize

\section*{Introduction}

Consider an $n$-dimensional real vector space with lattice of
integer points in it. The boundary of the convex hull of all
integer points contained inside one of the $n$-dimensional
invariant cones for a hyperbolic $n$-dimensional linear operator
without multiple eigenvalues is called a {\it sail} in the sense
of Klein. The set of all sails of such $n$-dimensional operator is
called {\it $(n{-}1)$-dimensional continued fraction} in the
sense of Klein (see in more details in Section~\ref{Ncfrac}). Any
sail is a polyhedral surface. In this work we study frequencies
of faces of multidimensional continued fractions.

There exists and is unique up to multiplication by a constant
function a form of the highest dimension on the manifold of
$n$-dimensional continued fractions in the sense of Klein, such
that the form is invariant under the natural action of the group
of projective transformations $PGL(n{+}1)$. A measure
corresponding to the integral of such form is called a {\it
M\"obius measure}. In the present paper we deduce an explicit
formulae to calculate invariant forms in special coordinates.
These formulae allow to give answers to some statistical
questions of theory of multi\-dimensional continued fractions. As
an example, we show in this work the results of approximate
calculations of frequencies for certain two-dimensional faces of
two-dimensional continued fractions.

\vspace{2mm}

A problem of generalization of ordinary continued fractions was
posed by C.~Hermite~\cite{Herm} in~1839. One of the most
interesting geometrical generalizations was introduced by F.~Klein
in~1895 in his works~\cite{Kle1} and~\cite{Kle2}. Unfortunately,
the computational complexity of multi\-di\-men\-sional continued
fractions did not allow to make significant advances in studies of
their properties one hundred years ago. V.~I.~Arnold originally
studying $A$-graded algebras~\cite{Arn1} faced with theory of
multidimensional continued fractions in the sense of Klein. Since
1989 he has formulated many problems on geometry and statistics of
multidimensional continued fractions, reviving an interest to the
study of multidimensional continued fractions (see the
works~\cite{ArnPT} and~\cite{Arn2}).

Multidimensional continued fractions in the sense of Klein are in
use in different branches of mathematics.
J.-O.~Moussafir~\cite{Mou1} and O.~N.~German~\cite{Ger} studied
the connection between the sails of multidimensional continued
fractions and Hilbert bases. In~\cite{Tsu} H.~Tsuchihashi
established the relationship between periodic multidimensional
continued fractions and multidimensional cusp singularities. This
relationship generalizes the classical relationship between
ordinary continued fractions and two-dimensional cusp
singularities known before. The combinatorial topological
multidimensional generalization of Lagrange theorem for ordinary
continued fractions was obtained by E.~I.~Korkina in~\cite{Kor1}
and the corresponding algebraic generalization by G.~Lachaud,
see~\cite{Lac}.

A large number of examples of multidimensional periodic continued
fraction were constructed by E.~Korkina
in~\cite{Kor0},~\cite{Kor2}, and~\cite{Kor3}, G.~Lachaud
in~\cite{Lac}, and~\cite{Lac2}, A.~D.~Bruno and V.~I.~Parusnikov
in~\cite{BP}, and~\cite{Par2}, and also by the author
in~\cite{Kar1} and~\cite{Kar2}. A portion of these two-dimensional
continued fractions is introduced at the web-site~\cite{site} by
K.~Briggs. A few examples of three-dimensional continued
fractions in four-dimensional space were constructed by the
author in~\cite{Kar3D}. The algorithms for constructing
multidimensional continued fractions are described in the works
of R.~Okazaki~\cite{Oka}, J.-O.~Moussafir~\cite{Mou2} and the
author~\cite{Kar4}.

For the first time the statement on statistics of numbers as
elements of ordinary continued fractions was formulated by
K.~F.~Gauss in~\cite{Gau}. This statement (see in the first
section) was proven further by R.~O.~Kuzmin~\cite{Kuz}, and
further was proven one more time by P.~L\'evy~\cite{Lev}. Further
investigations in this direction were made by E.~Wirsing
in~\cite{Wir}. (A basic notions of theory of ordinary continued
fractions is described in the books~\cite{Hin} by A.~Ya.~Hinchin
and~\cite{Arn2} by V.~I.~Arnold.) In 1989 V.~I.~Arnold
generalized statistical problems to the case of one-dimensional
and multidimensional continued fractions in the sense of Klein,
see in~\cite{ArnProb}, \cite{ArnPT}, and~\cite{ArnDyn}.

One-dimensional case was studied in details by M.~O.~Avdeeva and
B.~A.~Bykovskii in the works~\cite{Avd1} and~\cite{Avd2}. In
two-dimensional and multidimensional cases V.~I.~Arnold
formulated many problems on statistics of sail characteristics of
multidimensional continued fractions such as an amount of
triangular, quadrangular faces and so on, such as their integer
areas, and length of edges, etc. A major part of these problems
is open nowadays, while some are almost completely solved.

M.~L.~Kontsevich and Yu.~M.~Suhov in their work~\cite{Kon} proved
the existence of the mentioned above statistics. Recently
V.~A.~Bykovskii announced some solutions made by his students of
some problems of V.~I.~Arnold, including a problem on what is much
frequent: triangles or quadrangles. At present paper we write
down in special coordinates a natural M\"obius measure of the
manifold of all $n$-dimensional continued fractions in the sense
of Klein. In particular, this allows to make approximate
calculations of relative frequencies of multidimensional faces of
multidimensional continued fractions.

Note that the M\"obius measure is used also in theory of energies
of knots and graphs, see in the works of Freedman~M.~H.,
He~Z.~-H., and Wang~Z.~\cite{Freed}, J.~O'Hara~\cite{O-H4} and the
author~\cite{KarGr}. For the case of one-dimensional continued
fractions the M\"obius measure is induced by the relativistic
measure of three-dimensional de~Sitter world.

This work is organized as follows. In the first section we give
necessary notions of theory of ordinary continued fractions. In
particular, we give the definition of Gauss-Kuzmin statistics.
Further in the second section we describe the smooth manifold
structure for the set of all $n$-dimensional continued fractions
and define M\"obius measure on it. In the third section we study
relative frequencies of faces of one-dimensional continued
fractions. These frequencies are proportional to the frequencies
of Gauss-Kuzmin statistics. In the fourth section we study
relative frequencies of faces of multidimensional continued
fractions. Finally, in the fifth section we show approximate
calculation results of relative frequencies for some faces of
two-dimensional continued fractions.

The author is grateful to V.~I.~Arnold for constant attention to
this work and useful remarks, and Mathematisch Instituut of
Universiteit Leiden for the hospitality and excellent working
conditions.

\section{One-dimensional continued fractions and Gauss-Kuzmin statistics}

Let $\alpha$ be an arbitrary rational. Suppose that
$$
\alpha=a_0+1/(a_1+1/(a_2+\ldots+1/(a_{n-1}+1/a_n)\ldots)),
$$
where $a_0$ is integer, and the remaining $a_i$, $i=1,\ldots, n$
are positive integers. An expression on the right side of this
equality is called  a {\it decomposition of $\alpha$ into a finite
ordinary continued fraction} and denoted by
$[a_0,a_1,\ldots,a_n]$. If $n{+}1$ --- the total number of the
elements of the decomposition is even, then the continued
fraction is said to be {\it even}, and if this number is odd,
then the continued fraction is said to be {\it odd}.

Let $a_0$ be integer, and $a_1, \ldots, a_n, \ldots$ be infinite
sequence of positive integers. Denote by $r_n$ the rational
$[a_0,\ldots, a_{n-1}]$. For such integers $a_i$, the sequence
$(r_n)$ always converges to some real $\alpha$. The limit
$$
\lim\limits_{n\to \infty}[a_0,a_1,\ldots,a_{n-1}]
$$
is called  the {\it decomposition of $\alpha$ into a infinite
ordinary continued fraction} and denoted by
$[a_0,a_1,a_2,\ldots]$.

Ordinary continued fractions possess the following basic
properties.
\begin{proposition}
{\bf a).} Any rational has exactly two distinct decompositions
into a finite ordinary continued fraction, one of them is even,
and the other is odd.\\
{\bf b).} Any irrational  has a unique decomposition into an
infinite ordinary continued fraction.\\
{\bf c).} A decomposition into finite ordinary continued fraction is rational.\\
{\bf d).} A decomposition into infinite ordinary continued
fraction is irrational.
\end{proposition}

Notice, that for any finite continued fraction
$[a_0,a_1,\ldots,a_n]$, where $a_n{\ne} 1$, the following holds:
$$
[a_0,a_1,\ldots,a_n]=[a_0,a_1,\ldots,a_n{-}1,1].
$$
This equality determines a one-to-one correspondence between the
sets of even and odd finite continued fractions.

Let $\alpha$ be some irrational between zero and unity, and let
$[0,a_1,a_2,a_3,\ldots]$ be its ordinary continued fraction.
Denote by $z_n(\alpha)$ the real $[0,a_n,a_{n+1},a_{n+2},\ldots]$.

Let $m_n(\alpha)$ denote the measure of the set of reals $\alpha$
contained in the segment $[0;1]$, such that $z_n(\alpha)<x$. In
his letters to P.~S.~Laplace K.~F.~Gauss formulated without proofs
the following theorem. It was further proved by
R.~O.~Kuzmin~\cite{Kuz}, and then proved one more time by
P.~L\'evy~\cite{Lev}.
\begin{theorem}{\bf Gauss-Kuzmin.}
For $0\le x\le 1$ the following holds:
$$
\lim\limits_{n\to\infty} m_n(x)=\frac{\lg(1+x)}{\lg 2}.
$$
\end{theorem}

Denote by $P_n(k)$ for an arbitrary integer $k>0$ the measure of
the set of all reals $\alpha$ of the segment $[0;1]$, such that
each of them has the number $k$ at $n$-th position. A limit $
\lim\limits_{n\to\infty} P_n(k)$ is called a {\it frequency of}
$k$ for ordinary continued fractions and denoted by $P(k)$.

\begin{corollary}
For any positive integer $k$ the following holds
$$
P(k)=\frac{1}{\ln 2}\ln\left(1+\frac{1}{k(k+2)}\right).
$$
\end{corollary}

\begin{proof}
Notice, that $P_n(k)=m_n(\frac{1}{k})-m_n(\frac{1}{k+1})$. Now
the statement of the corollary follows from Gauss-Kuzmin theorem.
\end{proof}

The problem of V.~I.~Arnold on the asymptotic behaviours of
frequencies of integers as elements of ordinary continued
fractions for rationals with bounded numerators and denominators
was completely studied by V.~A.~Bykovskii and M.~O.~Avdeevain the
works~\cite{Avd1} and~\cite{Avd2}. It turns out that such
frequencies coincide with frequencies $P(k)$ defined above.

\section{Multidimensional continued fractions in the sense of Klein}\label{Ncfrac}

\subsection{Geometry of ordinary continued fractions}

Consider a two-dimensional plane with standard Euclidean
coordinates. A point is said to be {\it integer}, if both its
coordinates are integers. An {\it integer length} of the segment
$AB$ with integer vertices $A$ and $B$ is the ration of its
Euclidean length and the minimal Euclidean length for integer
vectors contained in the segment $AB$, we denote it
by~$\lil(AB)$. An {\it integer $($non-oriented$)$ area} of the
polygon $P$ is the ratio of its Euclidean area and the minimal
Euclidean area for the triangles with integer vertices, we denote
it by~$\ls(P)$. The quantity $\ls(P)$ coincides with doubled
Euclidean area of the polygon $P$.

\begin{figure}
$$\epsfbox{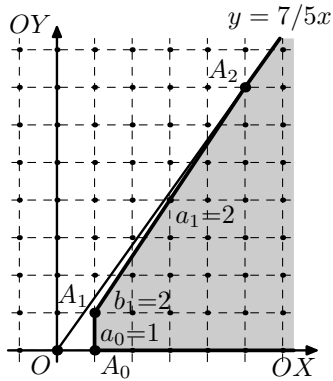}$$
\caption{The sail for the continued fraction of
$7/5=[1,2,2]$.}\label{7_5.1}
\end{figure}

For an arbitrary real $\alpha\ge 1$ we consider an angle in the
first orthant defined by the rays $\{(x,y)|y=0, x\ge 0\}$ and
$\{(x,y)|y=\alpha x, x\ge 0\}$. The boundary of the convex hull
of the set of all integer points in the closure of this angle
except the origin $O$ is a broken line, consisting of segments and
possible of a ray or two rays contained in the sides of the
angle. The union of all segments of that broken line is called
the {\it sail} of the angle. The sail of the angle is a finite
broken line for rational $\alpha$ and it is an infinite broken
line for irrationals. Denote the point with coordinates $(1,0)$
by $A_0$, and denote all the others vertices of the broken line
consequently by $A_1$, $A_2, \ldots$ Let $a_i=\lil(A_iA_{i+1})$
for $i=0,1,2,\ldots$, let also $b_i=\ls(A_{i-1}A_{i}A_{i+1})$ for
$i=1,2,3,\ldots$, then the following equality holds
$$
\alpha=[a_{0},b_{1}, a_{1}, b_{2}, a_{2}, b_{3}, a_{3}, \ldots ].
$$
On Figure~\ref{7_5.1} we examine an example of
$\alpha=7/5=[1,2,2]$.

\subsection{Definition of multidimensional continued fractions}\label{multCF}

Based on geometrical const\-ruction that we describe in the
previous subsection F.~Klein introduced the following geometrical
gene\-ralization of ordinary continued fractions to the
mul\-ti\-dimensional case (see~\cite{Kle1} and~\cite{Kle2}).

Consider arbitrary $n{+}1$ hyperplanes in~$\r^{n+1}$, such that
their intersection consists of a unique point --- of the origin.
The complement to the union of these hyperplanes consist of
$2^{n+1}$ open orthants. Consider one of them. The boundary of
the convex hull for the set of all integer points of the closure
of the orthant except the origin is called the {\it sail} of the
orthant. The set of all $2^{n+1}$ sails is called the {\it
$n$-dimensional continued fraction}, related to the given $n{+}1$
hyperplanes. An intersection of a hyperplane with the sail is
said to be a {\it $k$-dimensional face of the sail} if it is
contained in some $k$-dimensional plane and is homeomorphic to
$k$-dimensional disc. (See also~\cite{Kar2}.)

Two multidimensional faces of multidimensional continued
fractions are said to be {\it integer-linear $($-affine$)$
equivalent}, if there exist a linear (affine) integer lattice
preserving transformation taking one face to the other. A class
of all integer-linear (-affine) equivalent faces is called an
{\it integer-linear $($-affine$)$ type} of any face of this class.

Let us define one useful integer-linear invariant of a plane.
Consider an arbitrary $k$-dimensional plane $\pi$ not containing
the origin, whose integer vectors generates a sublattice of rank
$k$ in the lattice of all integer vectors. Let the Euclidean
distance from the origin to the plane $\pi$ equal $\ell$. Denote
by $\ell_0$ the minimal nonzero Euclidean distance to $\pi$ from
integer points of the plane (of dimension $k{+}1$) spanning the
given plane $\pi$ and the origin. The ratio $\ell/\ell_0$ is
called the {\it integer distance} from the origin to the
plane~$\pi$.

Let us now describe one of the original problems of V.~I.~Arnold
on statistics of faces of multidimensional continued fractions.
Note that for any real hyperbolic operator with distinct
eigenvalues there exists a unique corresponding multidimensional
continued fraction. One should take invariant hyperplanes for the
action of the operator as hyperplanes defining the corresponding
multidimensional continued fraction. Let us consider only
three-dimensional hyperbolic operators, that are defined by
integer matrices with rational eigenvalues. Denote the set of all
such operators by $A_3$. A continued fraction for any operator of
$A_3$ consists of finitely many faces. Denote by $A_3(m)$ the set
of all the operators of $A_3$ with bounded above by $m$ sums of
absolute values of all its coefficients. The number of such
operators is finite. Let us calculate the number of triangles,
quadrangles and so on among continued fractions, constructed for
the operators of $A_3(m)$. While $m$ tends to infinity we have a
general distribution of the frequencies for triangles, quadrangles
and so on. The problem of V.~I.~Arnold includes the study of the
properties of such distribution (for instant, {\it what is more
frequent: triangles or quadrangles}, {\it what is the frequency
of integer points inside the faces}, etc.). Note that this
problem still has not been completely studied. Surely, the
questions formulated above can be easily generalized to the
multidimensional case.

V.~I.~Arnold has also formulated statistical problems for special
algebraic periodic mul\-ti\-dimensional continued fractions. For
more information see~\cite{ArnPT} and~\cite{ArnDyn}.

\subsection{Smooth manifold of $n$-dimensional continued
fractions}

Denote the set of all continued fractions of dimension $n$ by
$CF_n$. Let us describe a natural structure of a smooth
nonsingular non-closed manifold on the set $CF_n$.

Consider an arbitrary continued fraction, that is defined by
unordered collection of hyperplanes $(\pi_1, \ldots, \pi_{n+1})$.
The enumeration of planes here is relative, without any ordering.
Denote by $l_i$ for $i=1,\ldots, n{+}1$ the intersection of all
the above hyperplanes except the hyperplane $\pi_i$. Obviously,
$l_1, \ldots, l_{n+1}$ are {\it independent} straight lines (i.e.
they are not contained in a hyperplane) passing through the
origin. These straight lines form an unordered collection of
independent straight lines. From the other side, any unordered
collection of $n{+}1$ independent straight lines uniquely
determines some continued fraction.

Denote the sets of all ordered collections of $n{+}1$ independent
and dependent straight lines by $FCF_n$ and $\Delta_n$
respectively. We say that $FCF_n$ is a space of $n$-dimensional
{\it framed continued fractions}. Also denote by $S_{n+1}$ the
permutation group acting on ordered collections of $n{+}1$
straight lines. In this notation we have:
$$
FCF_n=\bigl(\underbrace{\r P^{n}\times \r P^{n}\times \ldots
\times \r P^{n}}_{\hbox{$n{+}1$ times}}\bigr) \setminus \Delta_n
\qquad \hbox{and} \qquad CF_n=FCF_n\bigr/ S_{n+1}.
$$
Therefore, the sets $FCF_n$ and $CF_n$ admit natural structures
of smooth manifolds that are identified by the structure of the
Cartesian product of $n{+}1$ projective spaces $\r P^{n}$. Note
also, that $FCF_n$ is an $((n{+}1)!)$-fold covering of $CF_n$. We
call the map of ``forgetting'' of the order in the ordered
collections the {\it natural projection} of the manifold $FCF_n$
to the manifold $CF_n$ and denote it $p$, $p:FCF_n\to CF_n$.

\subsection{M\"obius measure on the manifolds of continued fractions}

A group $PGL(n{+}1,\r)$ of transformations of $\r P^n$ takes the
set of all straight lines passing through the origin of
$(n{+}1)$-dimensional space into itself. Hence, $PGL(n{+}1,\r)$
naturally acts on the manifolds $CF_n$ and $FCF_n$. Furthermore,
the action of $PGL(n{+}1,\r)$ is transitive, i.~e. it takes any
(framed) continued fraction to any other. Note that for any
$n$-dimensional (framed) continued fraction the subgroup of
$PGL(n{+}1,\r)$ taking this continued fraction to itself is of
dimension $n$.

\begin{definition}
A form of the manifold $CF_n$ (respectively $FCF_n$) is said to be
a {\it M\"obius form} if it is invariant under the action of
$PGL(n{+}1,\r)$.
\end{definition}

Transitivity of the action of $PGL(n{+}1,\r)$ implies that all
$n$-dimensional M\"obius forms of the manifolds $CF_n$ and
$FCF_n$ are proportional if exist.

Let $\omega$ be some volume form of the manifold $M$. Denote by
$\mu_\omega$ a measure of the manifold $M$ that  at any open
measurable set $S$ contained at the same piece-wise connected
component of $M$ is defined by an equality:
$$
\mu_\omega(S)=\left|\int_S \omega\right|.
$$

\begin{definition}
A measure $\mu$ of the manifold $CF_n$ ($FCF_n$) is said to be a
{\it M\"obius measure} if there exist a M\"obius form $\omega$ of
$CF_n$ ($FCF_n$) such that $\mu=\mu_\omega$.
\end{definition}

Note that any two M\"obius measures of $CF_n$ ($FCF_n$) are
proportional.

\begin{remark}
The projection $p$ projects the M\"obius measures of the manifold
$FCF_n$ to the M\"obius measures of the manifold $CF_n$. That
establishes an isomorphism between the spaces of M\"obius
measures for $CF_n$ and $FCF_n$. Since the manifold of framed
continued fractions possesses simpler chart system, all formulae
of the work are given for the case of framed continued fractions
manifold. To calculate a measure of some set $F$ of the unframed
continued fractions manifold one should: take $p^{-1}(F)$;
calculate M\"obius measure of the obtained set of the manifold of
framed continued fractions; divide the result by $(n{+}1)!$.
\end{remark}

\section{One-dimensional case}

\subsection{Explicit formulae for the M\"obius form}

Let us write down M\"obius forms of the framed one-dimensional
continued fractions manifold $FCF_1$ explicitly in special charts.

Consider a vector space $\r^2$ equipped with standard metrics on
it. Let $l$ be an arbitrary straight line in $\r^2$ that does not
pass through the origin, let us choose some Euclidean coordinates
$O_lX_l$ on it. Denote by $FCF_{1,l}$ a chart of the manifold
$FCF_1$ that consists of all ordered pairs of straight lines both
intersecting $l$. Let us associate to any point of $FCF_{1,l}$
(i.~e. to a collection of two straight lines) coordinates
$(x_l,y_l)$, where $x_l$ and $y_l$ are the coordinates on $l$ for
the intersections of $l$ with the first and the second straight
lines of the collection respectively. Denote by $|\bar{v}|_l$ the
Euclidean length of a vector $\bar{v}$ in the coordinates
$O_lX_lY_l$ of the chart $FCF_{1,l}$. Note that the chart
$FCF_{1,l}$ is a space $\r \times \r$ minus its diagonal.

Consider the following form in the chart $FCF_{1,l}$:
$$
\omega_l(x_l,y_l)=\frac{dx_l\wedge dy_l}{|x_l-y_l|_l^2}.
$$

\begin{proposition}
The measure $\mu_{\omega_l}$ coincides with the restriction of
ssome M\"obius measures to $FCF_{1,l}$.
\end{proposition}

\begin{proof}
Any transformation of the group $PGL(2,\r)$ is in the one-to-one
cor\-res\-pon\-dence with the set of all projective
transformations of the straight line~$l$ projectivization. Note
that the expression
$$
\frac{\Delta x_l\Delta y_l}{|x_l-y_l|_l^2}
$$
is an infinitesimal cross-ratio of four point with coordinates
$x_l$, $y_l$, $x_l{+}\Delta x_l$ and $y_l{+}\Delta y_l$. Hence
the form $\omega_l(x_l,y_l)$ is invariant for the action of
transformations (of the everywhere dense set) of the chart
$FCF_{1,l}$, that are induced by projective transformations of
$l$. Therefore, the measure $\mu_{\omega_l}$ coincides with the
restriction of some M\"obius measures to $FCF_{1,l}$.
\end{proof}

\begin{corollary}
A restriction of an arbitrary M\"obius measure to the chart
$FCF_{1,l}$ is proportional to $\mu_{\omega_l}$.
\end{corollary}

\begin{proof}
The statement follows from the proportionality of any two
M\"obius measures.
\end{proof}

Consider now the manifold $FCF_1$ as a set of ordered pairs of
distinct points on a circle $\r/\pi\z$ (this circle is a
one-dimensional projective space obtained from unit circle by
identifying antipodal points). The doubled angular coordinate
$\varphi$ of the circle $\r/\pi\z$ inducing by the coordinate $x$
of straight line $\r$ naturally defines the coordinates
$(\varphi_1, \varphi_2)$ of the manifold $FCF_1$.

\begin{proposition}\label{mob2}
The form $\omega_l(x_l,y_l)$ is extendable to some form $\omega_1$
of $FCF_1$. In coordinates $(\varphi_1, \varphi_2)$ the form
$\omega_1$ can be written as follows:
$$
\omega_1=\frac{1}{4}\cot
^2\left(\frac{\varphi_1-\varphi_2}{2}\right) d\varphi_1\wedge
d\varphi_2.
$$
\end{proposition}

We leave a proof of Proposition~\ref{mob2} as an exercise for the
reader.

\subsection{Relative frequencies of faces of one-dimensional continued fractions}

Without loose of generality in this subsection we consider only
M\"obius form $\omega_1$ of Proposition~\ref{mob2}. Denote the
natural projection of the form $\mu_{\omega_1}$ to the manifold of
one-dimensional continued fractions $CF_1$ by $\mu_1$.

Consider an arbitrary segment $F$ with vertices at integer
points. Denote by $CF_1(F)$ the set of continued fractions that
contain the segment $F$ as a face.

\begin{definition} The quantity $\mu_1 (CF_1(F))$ is called
{\it relative frequency} of the face $F$.
\end{definition}

Note that the relative frequencies of faces of the same
integer-linear type are equivalent. Any face of one-dimensional
continued fraction is at unit integer distance from the origin.
Thus, integer-linear type of a face is defined by its integer
length (the number of inner integer points plus unity). Denote
the relative frequency of the edge of integer length $k$ by
$\mu_1(''k'')$.

\begin{proposition}
For any positive integer $k$ the following holds:
$$
\mu_1(''k'')=\ln \left( 1+\frac{1}{k(k+2)}\right).
$$
\end{proposition}

\begin{proof}
Consider a particular representative of an integer-linear type of
the length $k$ segment: the segment with vertices $(0,1)$ and
$(k,1)$. One-dimensional continued fraction contains the segment
as a face iff one of the straight lines defining the fraction
intersects the interval with vertices $(-1,1)$ and $(0,1)$ while
the other straight line intersects the interval with vertices
$(k,1)$ and $(k{+}1,1)$, see on Figure~\ref{conf1}.

\begin{figure}
$$\epsfbox{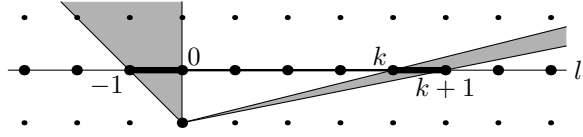}$$
\caption{Rays defining a continued fraction should lie in the
domain colored in gray.}\label{conf1}
\end{figure}

For the straight line $l$ defined by the equation $y=1$ we
calculate the M\"obius measure of Cartesian product of the
described couple of intervals. By the last subsection it follows
that this quantity coincides with relative frequency
$\mu_1(''k'')$. So,
$$
\begin{array}{l}
\displaystyle \mu_1(''k'')=\int\limits_{-1}^{0}
\int\limits_{k}^{k+1} \frac{dx_l dy_l}{(x_l-y_l)^2}=
\int\limits_{k}^{k+1} \left(
\frac{1}{y_l} -\frac{1}{y_l+1} \right) dy_l =\\
\displaystyle
 \ln \left( \frac{(k+1)(k+1)}{k(k+2)}\right)=
 \ln \left(1+\frac{1}{k(k+2)}\right).
\end{array}
$$
This proves the proposition.
\end{proof}

\begin{remark} Note that the argument of the logarithm
$\frac{(k+1)(k+1)}{k(k+2)}$ is a cross-ratio of points $(-1,1)$,
$(0,1)$, $(k,1)$, and $(k{+}1,1)$.
\end{remark}

\begin{corollary}
Relative frequency $\mu_1(''k'')$ up to the factor
$$ \ln 2=\int\limits_{-1}^{0}\int\limits_{1}^{+\infty}\frac{dx_l
dy_l}{(x_l-y_l)^2}
$$
coincides with Gauss-Kuzmin frequency $P(k)$ for $k$ to be an
element of continued fraction. \qed
\end{corollary}

\section{Multidimensional case}

\subsection{Explicit formulae for the M\"obius form}

Let us now write down explicitly M\"obius forms for the manifold
of framed $n$-dimensional continued fractions $FCF_n$ for
arbitrary $n$.

Consider $\r^{n+1}$ with standard metrics on it. Let $\pi$ be an
arbitrary hyperplane of the space $\r^{n+1}$ with chosen
Euclidean coordinates $OX_1\ldots X_{n}$, let also $\pi$ does not
pass through the origin. By the chart $FCF_{n,\pi}$ of the
manifold $FCF_n$ we denote the set of all collections of $n{+}1$
ordered straight lines such that any of them intersects $\pi$.
Let the intersection of $\pi$ with $i$-th plane is a point with
coordinates $(x_{1,i},\ldots,x_{n,i})$ at the plane $\pi$. For an
arbitrary tetrahedron $A_1\dots A_{n+1}$ in the plane $\pi$ we
denote by $V_\pi(A_1\dots A_{n+1})$ its oriented Euclidean volume
in the coordinates $OX_{1,1}\ldots X_{n,1}X_{1,2}\ldots
X_{n,n+1}$ of the chart $FCF_{n,\pi}$. Denote by $|\bar{v}|_\pi$
the Euclidean length of the vector $\bar{v}$ in the coordinates
$OX_{1,1}\ldots X_{n,1}X_{1,2}\ldots X_{n,n+1}$ of the chart
$FCF_{n,\pi}$. Note that the map $FCF_{n,\pi}$ is everywhere
dense in $(\r^{n})^{n+1}$.

Consider the following form in the chart $FCF_{n,\pi}$:
$$
\omega_\pi(x_{1,1},\ldots,x_{n,n+1})=\frac{
\bigwedge\limits_{i=1}^{n+1} \left(\bigwedge\limits_{j=1}^{n}
dx_{j,i}\right) }
{V_\pi(A_1\dots A_{n+1})^{n+1}}.
$$

\begin{proposition}
The measure $\mu_{\omega_\pi}$ coincides with the restriction of
some of M\"obius measure to $FCF_{n,\pi}$.
\end{proposition}

\begin{proof}
Any transformation of the group $PGL(n{+}1,\r)$ is in the
one-to-one cor\-res\-pon\-dence with the set of all projective
transformations of the plane~$\pi$. Let us show that the form
$\omega_\pi$ is invariant for the action of transformations (of
the everywhere dense set) of the chart $FCF_{n,\pi}$, that are
induced by projective transformations of hyperplane $\pi$.

Let us at each point of the tangent space to $FCF_{n,\pi}$ define
a new basis corresponding to the directions of edges of the
corresponding tetrahedron in $\pi$. Namely, consider an arbitrary
point $(x_{1,1},\ldots,x_{n,n{+}1})$ of the chart $FCF_{n,\pi}$
and the tetrahedron $A_1\dots A_{n+1}$ in hyperplane $\pi$
corresponding to the point. Let
$$
\bar{f}_{ij}=\frac{\bar{A_jA_i}}{|\bar{A_jA_i}|_\pi}, \qquad
i,j=1,\ldots,n{+}1; \quad i\ne j.
$$

The basis constructed above continuously depends on the point of
the chart $FCF_{n,\pi}$. By $dv_{ij}$ we denote the 1-form
corresponding to the coordinate along the vector $\bar{f}_{ij}$
of $FCF_{n,\pi}$.

Denote by $A_i=A_i(x_{1,i},\ldots,x_{n,i})$ the point depending
on the coordinates of the plane $\pi$ with coordinates
$(x_{1,i},\ldots,x_{n,i})$, $i=1,\ldots,n{+}1$. Let us rewrite
the form $\omega_\pi$ in new coordinates.
$$
\begin{array}{c}
\omega_\pi(x_{1,1},\ldots,x_{n,n{+}1})=\\
\displaystyle \prod\limits_{i=1}^{n+1}\left(
\frac{V_{\pi}(A_i,A_1,\ldots,A_{i-1},A_{i+1},\ldots A_{n+1})}
{\prod\limits_{k=1,k\ne i}^{n+1} |\bar{A_kA_i}|_\pi}\right)
 \cdot\frac{dv_{21}\wedge dv_{31}\wedge \ldots \wedge
dv_{n,n{+}1}}{V_\pi(A_1\ldots A_{n+1})^{n+1}}=\\
\displaystyle (-1)^{[\frac{n+3}{4}]}\cdot
\frac{dv_{21}\wedge dv_{12}}{|\bar{A_1A_2}|_\pi^2} \wedge
\frac{dv_{32}\wedge dv_{23}}{|\bar{A_2A_3}|_\pi^2} \wedge \ldots
\frac{dv_{n{+}1,n}\wedge dv_{n,n{+}1}}{|\bar{A_{n+1}A_n}|_\pi^2},
\end{array}
$$
here by $[a]$  we denote the maximal integer not exceeding $a$.

Like in one-dimensional case the expression
$$
\frac{\Delta v_{ij}\Delta v_{ji}}{|\bar{A_iA_j}|^2}
$$
for the infinitesimal small $\Delta v_{ij}$ and $\Delta v_{ji}$
is the infinitesimal cross-ratio of four points: $A_i$, $A_j$,
$A_i{+}\Delta v_{ji}\bar{f}_{ji}$, and $A_j{+}\Delta v_{ij}
\bar{f}_{ij}$ of the straight line $A_iA_j$. Therefore, the form
$\omega_\pi$ is invariant for the action of transformations (of
the everywhere dense set) of the chart $FCF_{n,\pi}$, that are
induced by projective transformations of hyperplane $\pi$. Hence
the measure $\mu_{\omega_\pi}$ coincides with the restriction of
some M\"obius measure to $FCF_{n,\pi}$.
\end{proof}

\begin{corollary}
A restriction of an arbitrary M\"obius measure to the chart
$FCF_{n,\pi}$ is proportional to $\mu_{\omega_\pi}$.
\end{corollary}

\begin{proof}
The statement follows from the proportionality of any two
M\"obius measures.
\end{proof}

Let us fix an origin $O_{ij}$ for the straight line $A_iA_j$. The
integral of the form $dv_{ij}$ (respectively $dv_{ij}$) for the
segment $O_{ij}P$ defines the {\it coordinate} $v_{ij}$
($v_{ij}$) of the point $P$ contained in the straight line
$A_iA_j$. As in one-dimensional case consider a projectivization
of the straight line $A_iA_j$. Denote the angular coordinates by
$\varphi_{ij}$ and $\varphi_{ji}$ respectively. In this
coordinates it holds:
$$
\frac{dv_{ij}\wedge dv_{ji}}{|\bar{A_iA_j}|_\pi^2}=
\frac{1}{4}\cot ^2\left(\frac{\varphi_{ji}-\varphi_{ij}}{2}\right)
d\varphi_{ij}\wedge d\varphi_{ji}.
$$
Then, the following is true.
\begin{corollary}\label{mobn}
The form $\omega_\pi$ extends to some form $\omega_n$ of $FCF_n$.
In coordinates $v_{ij}$ the form $\omega_n$ is as follows:
$$
\omega_n=\frac{(-1)^{[\frac{n+3}{4}]}}{2^{n(n+1)}}\left(\prod\limits_{i=1}^{n+1}\prod\limits_{j=i+1}^{n+1}
\cot
^2\left(\frac{\varphi_{ij}-\varphi_{ji}}{2}\right)\right)\cdot
\left(\bigwedge\limits_{i=1}^{n+1}\left(\bigwedge\limits_{j=i+1}^{n+1}
d\varphi_{ij}\wedge d\varphi_{ji}\right) \right).
$$
\end{corollary}

\subsection{Relative frequencies of faces of multidimensional continued fractions}

As in one-dimensional case without loose of generality we
consider the form $\omega_n$ of Corollary~\ref{mobn}. Denote by
$\mu_n$ the projection of the measure $\mu_{\omega_n}$ to the
manifold of multidimensional continued fractions $CF_n$.

Consider an arbitrary polytope $F$ with vertices at integer
points. Denote by $CF_n(F)$ the set of $n$-dimensional continued
fractions that contain the polytope $F$ as a face.

\begin{definition} The value $\mu (CF_n(F))$ is called
the {\it relative frequency} of a face $F$.
\end{definition}

Relative frequencies of faces of the same integer-linear type are
equivalent.

\begin{problem}\label{p1}
Find integer-linear types of $n$-dimensional continued fractions
with the highest relative frequencies. Is it true that the number
of integer-linear types of faces with relative frequencies bounded
above by some constant is finite? Find its asymptotics for the
constant tending to infinity.
\end{problem}

Problem~\ref{p1} is open for $n\ge 2$.

\begin{conjecture}
Relative frequencies of faces are proportional to the frequencies
of faces in the sense of Arnold $($see in
Subsection~\ref{multCF}$)$.
\end{conjecture}
This conjecture is checked in the present work for the case of
one-dimensional continued fractions. It is still open for the
$n$-dimensional case for $n\ge 2$.

\section{Examples of calculation of relative frequencies for faces
in two-dimensional case}

\subsection{A method of relative frequencies computation}

Let us describe a method of relative frequencies computation in
two-dimensional case more detailed.

Consider a space $\r^3$ with standard metrics on it. Let $\pi$ be
an arbitrary plane in $\r^3$ not passing through the origin and
with fixed system of Euclidean coordinates $O_\pi X_\pi Y_\pi$.
Let $FCF_{2,\pi}$ be the corresponding chart of the manifold
$FCF_2$ (see the previous section). For an arbitrary triangle
$ABC$ of the plane $\pi$ we denote by $S_\pi(ABC)$ its oriented
Euclidean area in the coordinates $O_\pi X_1Y_1X_2Y_2X_3Y_3$ of
the chart $FCF_{2,\pi}$. Denote by $|\bar{v}|_\pi$ the Euclidean
length of the vector $\bar{v}$ in the coordinates $O_\pi
X_1Y_1X_2Y_2X_3Y_3$ of the chart $FCF_{2,\pi}$. Consider the
following form in the chart $FCF_{2,\pi}$:
$$
\omega_\pi(x_1,y_1,x_2,y_2,x_3,y_3)=\frac{dx_1\wedge dy_1\wedge
dx_2\wedge dy_2 \wedge  dx_3\wedge
dy_3}{S_\pi((x_1,y_1)(x_2,y_2)(x_3,y_3))^3}.
$$

Note that the oriented area $S_\pi$ of the triangle
$(x_1,y_1)(x_2,y_2)(x_3,y_3)$ can be expressed in the coordinates
$x_i$, $y_i$ as follows:
$$
S_\pi((x_1,y_1)(x_2,y_2)(x_3,y_3))=\frac{1}{2}
\bigl(x_3y_2-x_2y_3+x_1y_3-x_3y_1+x_2y_1-x_1y_2\bigr).
$$

For the approximate computations of relative frequencies of faces
it is useful to rewrite the form $\omega_\pi$ in the dual
coordinates (see Remark~\ref{polezno}). Define a triangle $ABC$
in the plane $\pi$ by three straight lines $l_1$, $l_2$, and
$l_3$, where $l_1$ passes through $B$ and $C$, $l_2$  passes
through $A$, and $C$, and $l_3$  passes through $A$, and $B$.
Define the straight line $l_i$ ($i=1,2,3$) in $\pi$ by the
equation (preliminary we make a translation of $\pi$ in such a way
that the origin is taken to some inner point of the triangle)
$$
a_ix+b_iy=1
$$
in $x$ and $y$ variables. Then if we know the 6-tuple of numbers
$(a_1,b_1,a_2,b_2,a_3,b_3)$ we can restore the triangle in the
unique way.

\begin{proposition}\label{goodkoords}
In coordinates  $a_1, b_1, a_2, b_2, a_3, b_3$ the form
$\omega_\pi$ can be written as follows:
$$
-\frac{8 da_1\wedge d b_1 \wedge da_2\wedge d b_2 \wedge da_3
\wedge d b_3}{(a_3b_2-a_2b_3+a_1b_3-a_3b_1+a_2b_1-a_1b_2)^3}.
$$ \qed
\end{proposition}

So, we reduce the computation of relative frequency for the face
$F$, i.~e. the value of $\mu_2(CF_2(F))$ to the computation of
measure $\mu_{\omega_2}(p^{-1}(CF_2(F)))$. Consider some plane
$\pi$ in $\r^3$ not passing through the origin. By
Corollary~\ref{mobn}
$$
\mu_{\omega_2}(p^{-1}(CF_2(F)))=\mu_{\omega_\pi}(p^{-1}(CF_2(F))\cap
(FCF_{2,\pi}).
$$
Finally the computation should be made for the set
$\mu_{\omega_\pi}(p^{-1}(CF_2(F))\cap (FCF_{2,\pi})$ in dual
coordinates $a_i, b_i$ (see Proposition~\ref{goodkoords}).

\begin{remark}\label{polezno} In $a_i, b_i$ coordinates the
computation of value of the relative frequency often reduces to
the estimation of the integral on the disjoint union of the
finite number of six-dimensional Cartesian products of three
triangles in $a_i, b_i$ coordinates (see
Proposition~\ref{goodkoords}). The integration over such a simple
domain greatly fastens the speed of approximate computations. In
particular, the integration can be reduced to the integration
over some 4-dimensional domain.
\end{remark}

\subsection{Some results}

In conclusion of the work we give some results of relative
frequencies calculation for some two-dimensional faces of
two-dimensional continued fractions.

Explicit calculations of relative frequencies for the faces seems
not to be realizable. Nevertheless it is possible to make
approximations of the corresponding integrals. Normally, the
greater area of the integer-linear type of the polygon is, the
lesser its relative frequency. The most complicated approximation
calculations correspond to the most simple faces, such as an
empty triangle.

\begin{figure}[ht]
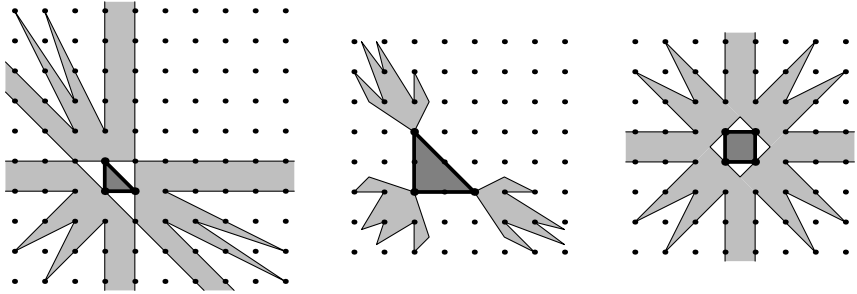

$$
\begin{array}{c}
\epsfbox{faces.1}
\end{array} \quad
\begin{array}{c}
\epsfbox{faces.3}
\end{array} \quad
\begin{array}{c}
\epsfbox{faces.2}
\end{array}
$$
\caption{The points painted in light-gray correspond to the
points at which the rays defining the two-dimensional continued
fraction can intersect the plane of the chosen face.}\label{faces}
\end{figure}

On Figure~\ref{faces} we show examples of the following faces:
triangular $(0,0,1)$, $(0,1,1)$, $(1,0,1)$ and $(0,0,1)$,
$(0,2,1)$, $(2,0,1)$ and quadrangular $(0,0,1)$, $(0,1,1)$,
$(1,1,1)$, $(1,0,1)$. For each face it is shown the plane
containing the face. The points painted in light-gray correspond
to the points at which the rays defining the two-dimensional
continued fraction can intersect the plane of the chosen face.

Faces of two-dimensional continued fractions for the majority of
integer-linear types lie at unit integer distance from the
origin. Only three infinite series and three partial examples of
faces lie at integer distances greater or equal to two from the
origin, see a detailed description in~\cite{KarPyr}
and~\cite{KarPyrNote}. If the distance to the face is increasing,
then the frequency of faces is reducing on average. The average
rate of reducing the frequency is unknown to the author.

\begin{table}[ht]
\begin{tabular}{|c|l|c|c|c||c|l|c|c|c|}
\hline {\bf N}$^\circ$ &
{\bf face}
& $\ls$ & $\ld$ &
{ $\mu_2$}
&
{\bf N}$^\circ$ &
{\bf face}
& $\ls$ & $\ld$ &
{ $\mu_2$}
\\
\hline\hline
{\bf I}$_1$ & $\begin{array}{c} \epsfbox{example.1} \end{array}$ &
3 & 1 & $1.3990\cdot 10^{-2}$
& {\bf VI}$_1$ & $\begin{array}{c} \epsfbox{example.6}
\end{array}$ & 7 & 1 & $ 3.1558\cdot 10^{-4}$
\\
\hline
{\bf I}$_3$ & $\begin{array}{c} \epsfbox{example.1} \end{array}$ &
3 & 3 & $1.0923\cdot 10^{-3}$
& {\bf VI}$_2$ & $\begin{array}{c} \epsfbox{example.6}
\end{array}$ & 7 & 2 & $ 3.1558\cdot 10^{-4}$
\\
\hline
{\bf II}$_1$ & $\begin{array}{c} \epsfbox{example.2} \end{array}$
& 5 & 1 & $1.5001\cdot 10^{-3}$
& {\bf VII}$_1$ & $\begin{array}{c} \epsfbox{example.7}
\end{array}$ & 11 & 1 & $ 3.4440\cdot 10^{-5}$
\\
\hline
{\bf III}$_1$ & $\begin{array}{c} \epsfbox{example.3} \end{array}$
& 7 & 1 & $3.0782\cdot 10^{-4}$
& {\bf VIII}$_1$ & $\begin{array}{c} \epsfbox{example.8}
\end{array}$ & 7 & 1 & $5.6828\cdot 10^{-4}$
\\
\hline
{\bf IV}$_1$ & $\begin{array}{c} \epsfbox{example.4} \end{array}$
& 9 & 1 &  $9.4173\cdot 10^{-5}$
& {\bf IX}$_1$ &  $\begin{array}{c} \epsfbox{example.9}
\end{array}$ & 7 & 1 & $1.1865\cdot 10^{-3}$
\\
\hline
{\bf V}$_1$ & $\begin{array}{c} \epsfbox{example.5} \end{array}$ &
11 & 1 & $3.6391\cdot 10^{-5}$
& {\bf X}$_1$ & $\begin{array}{c} \epsfbox{example.10}
\end{array}$ & 6 & 1 & $9.9275\cdot 10^{-4}$
\\
\hline
\end{tabular}
\caption{Some results of calculations of relative frequencies.
}\label{megatable}
\end{table}

In Table~\ref{megatable} we show the results of relative
frequencies calculations for 12 integer-linear types of faces. In
a column ``N$^\circ$'' we write a special sign for integer-affine
type of a face. The index denotes the integer distance from the
corresponding face to the origin. In a column ``face'' we draw a
picture of integer-affine type of the face. Further in a column
``$\ls$'' we write down integer areas of faces, and in a column
``$\ld$'' we write down integer distances from the planes of faces
to the origin. Finally in a column ``$\mu_2$'' we show the
results of the approximate relative frequency calculations for the
corresponding integer-linear types of faces.

Note that in the given examples the integer-affine type and
integer distance to the origin determines the integer-linear type
of the face.

In conclusion of this section we give two simple statements on
relative frequencies of faces.

\begin{statement}\label{1&2}
Faces of the same affine-linear type at integer distance to the
origin equivalent to $1$ and at integer distance to the origin
equivalent to $2$ always have the same relative frequencies $($see
for example {\bf VI$_1$} and {\bf VI$_2$} of
Table~\ref{megatable}$)$.
\end{statement}

Denote by $A_n$ the triangle with vertices $(0,0,1)$, $(n,0,1)$,
and $(0,n,1)$. Denote by $B_n$ the square with vertices
$(0,0,1)$, $(n,0,1)$, $(n,n,1)$, and $(0,n,1)$.

\begin{statement}
The following holds
$$
\lim\limits_{n\to \infty} \frac{\mu(CF_n(A_n))}{\mu(CF_n(B_n))}=8.
$$
\end{statement}

\end{document}